\documentclass{amsart}
\usepackage{epsf} 
\usepackage{amssymb,latexsym, amsmath, amscd, array} 
 
\swapnumbers 

\theoremstyle{plain} 
\newtheorem{thm}{Theorem}
\newtheorem{cory}[thm]{Corollary} 
\newtheorem{cor}[thm]{Corollary} 
\newtheorem{lem}[thm]{Lemma} 
\newtheorem{theorem}[thm]{Theorem} 
 
\newtheorem{prop}[thm]{Proposition} 

\newcommand\theoref{Theorem~\ref} 
\newcommand\lemref{Lemma~\ref} 
\newcommand\propref{Proposition~\ref} 
\newcommand\corref{Corollary~\ref}

\theoremstyle{definition}

\def\p{{\noindent \it Proof. }}

\def\ga{\alpha}

\def\Q{{\mathbb Q}} 
\def\R{{\mathbb R}}

\def\ra{{\R A}}

\def\pa{\partial}

\def\la{\langle}
\def\ra{\rangle}

\def\m{\medskip}

\begin{document} 

\title[Non-formal closed $(k-1)$-connected manifolds] {Examples of non-formal 
closed $(k-1)$-connected manifolds of dimensions $4k-1$ and more}

\author[A.~Dranishnikov and Yu.~Rudyak]
{Alex N. Dranishnikov and Yuli B. Rudyak} 
\address{A. Dranishnikov, Department of Mathematics, University 
of Florida, 358 
Little Hall, Gainesville, FL 32611-8105, USA} 
\email{dranish@math.ufl.edu} 
\address{Yu. Rudyak, Department of Mathematics, University 
of Florida, 358 
Little Hall, Gainesville, FL 32611-8105, USA} 
\email{rudyak@math.ufl.edu} 

\begin{abstract} We construct closed $(k-1)$-connected 
manifolds of dimensions $\ge 4k-1$ that possess non-trivial rational Massey 
triple products. We also construct examples of manifolds $M$ such that all the 
cup-products of elements of $H^k(M)$ vanish, while the group $H^{3k-1}(M;\Q)$ 
is generated by Massey products: such examples are useful for theory of systols.

\end{abstract}

\maketitle 

For every $k$ we construct closed $(k-1)$-connected manifolds of dimensions 
$\ge 4k-1$ that possess non-trivial rational Massey triple products and 
therefore are non-formal.  
For $k=1$ such manifolds can be obtained as the products of Heisenberg manifold 
with circles. For $k=2$ such examples are also known, see e.g. \cite{Oprea, 
FernandezMunoz}, but even in this case our construction seems more direct and 
simple. 

\m  Miller~\cite{Miller} proved that every 
closed $(k-1)$-connected manifold $M$ of dimension $\leq 4k-2$ is formal. In 
particular, all rational Massey products in $M$ vanish. So, neither Miller's 
nor our results can be improved.

\m                                                               
Given a diagram 
$$
\CD
B\supset A @>f>> Y
\endCD
$$
we denote by $Z_f$ its double mapping cylinder.

\m Recall that a subset $S$ of a space $\R^m$ is called {\it radial} if, for 
all points $s\in S$, the linear segment $[0,s]$ contains precisely one point of 
$S$ (namely, $s$). 
 
\begin{prop}\label{embed}
Let $B$ be a finite polyhedron in $\R^m, m>1$, let $A$ be a subpolyhedron of 
$B$ such that $A\setminus \{0\}$ is radial in $\R^m$, and let $Y$ be a finite 
polyhedron in $\R^n$. Then the double cylinder $Z_f$ of any 
simplicial map $f: A \to Y$ admits a PL embedding in $\R^{m+n}$. 
\end{prop}

\p We denote by $0_m$ and $0_n$ the origins of spaces $R^m$ and $R^n$, 
respectively. We first consider the case when $0_m\notin A$, We assume that $Y$ 
is far away from $0_n$. Let $\Gamma\in \R^{m}\times \R^n$ be the graph of the 
map $f$. We join every point $(x,f(x))\in \Gamma,x\in A$ with the point 
$(0_m,f(x))\in \R^m\times Y\subset \R^{m+n}$ by the linear segment. Then, since 
$A$ is radial, we get an embedding of the mapping cylinder $M_f$ of $f$ to 
$\R^{m+n}$. 
Moreover, if we join the points $(x,0_n)$ with $(x,f(x))$ by the linear segment, 
we still have an embedding $M_f\hookrightarrow \R^{m+n}$. Here (the image of) 
$M_f$ is formed by segments $[(x,0_n),(x,f(x))]$ and $[(x,f(x)), (0_m,f(x))]$. 
Finally, we get an embedding of the double mapping cylinder $Z_f$ to $\R^{m+n}$ 
by adding the space $B$ to the embedded mapping cylinder $M_f$.

\m The case $0_m\in A$ can be considered similarly. We can assume that there is 
a point $y_0\in Y$ which is the closest to $0_n\in \R^n$, i.e. $||y_0||< ||y||$ 
if $y\ne y_0$ and $y\in Y$. We can also assume that $f(0_m)= y_0$. Consider the 
map $f'=f|(A\setminus \{0\})$ and the embedding $i:Z_{f'} \to \R^{m+n}$ as 
above. 
Then $i(Z_{f'})\cup [0_m, y_0]$ is an embedding of $Z_f$. 
\qed

\begin{cory} Let $Y$ be a finite polyhedron in $\R^n$, and let $f: \vee_i 
S^{m-1}_i \to Y, i=1, \ldots, k$ be a simplicial map, where $S^{m-1}_i$ is the 
copy of the sphere $S^{m-1}$. Then the cone $C_f$ of $f$ can be simplicially 
embedded in $\R^{m+n}$.
\end{cory}

\p Choose a base point on the boundary of each disc $D^m_i, i=1, \ldots, k$
and consider the wedge $\vee_{i=1}^mD^m_i$. We can regard this wedge as a 
polyhedron in $\R^m$ such that the base point is the origin and $\vee 
S^{m-1}_i\setminus \{0\}$ is a radial set. Now the claim follows from 
\propref{embed}.
\qed

\m Consider the wedge $K=S^{k_1}\vee S^{k_2}\vee S^{k_3}$ of spheres with 
$k_i\ge 2$ and let 
$\iota_r\in \pi_{k_r}(K)$ be represented by the inclusion map 
$S^{k_r}\subset K$. Set $m=k_1+k_2+k_3-1$, let
$f: S^{m-1} \to K$ represent the element 
$[\iota_1,[\iota_2,\iota_3]]$, and let $X$ be the cone of the map $f$. Let 
$\ga_i\in H^{k_i}(X)$ be the cohomology class which takes the value 1 on the 
cell $S^k_i$ of $X$ and 0 on other cells. We recall the following classical 
result

\begin{theorem}\label{um}
The Massey product $\la \ga_1, \ga_2, \ga_3\ra 
\in H^{k_1+k_2+k_3-1}(X)$ has the zero indeterminacy and takes the value 
$(-1)^{k_1}$ on the $(m-1)$-dimensional cell of $X$. 
\end{theorem}

\p See \cite[Lemma 7]{UeharaMassey}.
\qed
 
\m Now let $k_1=k_2=k_3=k$ and consider the corresponding space $X$.
According to \propref{embed}, $X$ admits a PL embedding in $\R^N$ with $N\ge 
4k$.  Fix such an embedding and let $W$ be a closed regular neighborhood of $X$ 
in $\R^N$. So, $W$ is a manifold  with the boundary $V=\partial W$. Furthermore, 
$W$ has the homotopy type of $X$. (Notice that $W$ is a PL manifold by the 
construction, but without loss of generality we can assume that $W$ is smooth.)

\begin{prop} The manifold $V$ is $(k-1)$-connected.
\end{prop}

\p Consier a sphere $S^i, i<k$ in $V$. Since $W$ is $(k-1)$-connected, there 
exists a disk $D^{i+1}$ in $W$ with $\pa D^{i+1}=S^i$. Since $i+1+\dim X\le 
4k-1< N$, we can assume that $D^{i+1}\cap X =\emptyset$. But $V$ is a retract of 
$W\setminus X$, and thus $S^i$ bounds a disk in $V$.
\qed  

\begin{prop}\label{homology}
$H^i(W, V)=H_{N-i}(X)$.
\end{prop} 

\p We have
$$
H^i(W, V)=H_{N-i}(W)=H_{N-i}(X)
$$
where the first equality holds by the Poincar\'e duality, see e.g. 
Dold~\cite{Dold}.
\qed

\m Consider the map
$$
\CD
g: V@>i>> W @>r >> X
\endCD
$$
where $i$ is the inclusion and $r$ is a deformation retraction.
 
\begin{thm}
If $N\ne 5k-1,\, 6k-2$, then the Massey 
product $\la g^*\ga_1, g_*\ga_2, g^*\ga_3\ra$ has zero indeterminacy and is 
non-zero.
\end{thm}

\p Notice that $H_i(X)=0$ for $i\ne 0, k, 3k-1$. We have 
$H^{2k-1}(W)=H^{2k-1}(X)=0$ and $H^{2k}(W,V)=H_{n-2k}(X)=0$. Now, in view of the 
exactness of the sequence $H^{2k-1}(W) \to H^{2k-1}(V) \to H^{2k}(W,V)$ we have 
$H^{2k-1}(V)=0$, 
and therefore the indeterminacy of the Massey product is zero. Furthermore, the 
map $i^*:H^{3k-1}(W) \to H^{3k-1}(V)$ is injective since 
$H^{3k-1}(W,V)=H_{n-3k+1}(X)=0$. Thus, 
the map $g^*: H^{3k-1}(X) \to H^{3k-1}(V)$ is injective. But $g^*\la\ga_1, 
g_*\ga_2, 
g^*\ga_3\ra = \la g^*\ga_1, g_*\ga_2, g^*\ga_3\ra$ because both parts of the 
equality have zero indeteminaces.
\qed

\m Thus, we have examples of $(k-1)$-connected manifolds with non-trivial triple 
Massey product of dimensions $d\ge 4k-1$ but $d\ne 5k-2,\, 6k-3$. In order to 
construct an example in exceptional dimensions, just take the double of the 
manifold $W$ (or multiple by the sphere of the correspondent dimension if $k\ne 
2$).

\m When we put the first version of the paper into the e-archive, Mikhail Katz 
asked us if we can construct a closed manifold $M$ such that all the 
cup-products of elements of $H^k(M)$ vanish, while the group $H^{3k-1}(M;\Q)$ 
is generated by Massey products. Now we present such an example. 

\begin{lem}\label{wedge}
 Consider a wedge $X\vee Y$ and three elements $u,v,w\in H^*(X)$ such that 
$uv=0$, $u|Y=0=v|Y$ and $w|X=0$. Then all the Massey products $\la u,v,w\ra$, 
$\la u,w,v\ra$ and $\la w,u,v\ra$ are trivial, i.e. they contain the zero 
element.
\end{lem}

\p This follows from the following fact: If $A\in C^*(X\vee Y)$ and $B\in 
C^*(X\vee Y)$ are cochains with the supports in $X$ and $Y$, respectively, than 
their product is equal to zero. We leave the details to the reader.
\qed

\m Consider the wedge $S^k_1\vee S^k_2\vee S^k_3\vee S^k_4$ of $k$-dimensional 
spheres, $k>1$. Let $\iota_m\in \pi_k(S^k_m)$ be the generator. Set 
\begin{equation}\label{z}
Z=\left(\vee_{i=1}^4S^k_i\right)\cup_{f_1}e^{3k-1}
\end{equation}
 where $f_1: S^{3k-2}\to \vee_{i=1}^4S^k_i$ represents the homotopy class 
$[\iota_1,[\iota_2, \iota_3]]$. Let 
$\ga_i\in H^k(Z)$ be the cohomology class which takes the value 1 on the cell 
$S^k_i$ of $Z$ and 0 on other cells.

\begin{cor}\label{zero} 
If at least one of the indices $i,j,k$ is equal to $4$, then $\la 
\ga_i, \ga_j, \ga_k\ra=0$ in $X$.
\end{cor}

\p This follows directly from \lemref{wedge} since 
$$
Z=\left(\left(\vee_{i=1}^3S^k_i\right)\cup_{f_1}e^{3k-1}\right)\vee S^k_4.
$$
\qed

\m For convenience of notation, we set $\iota_5=\iota_1$ and $\iota_6=\iota_2$ . 
Let $f_m: 
S^{3k-2} \to \vee_{i=1}^4S^k_i$ be the map which represents $[\iota_{m}, 
[\iota_{m+1}, \iota_{m+2}]], m=1,2,3,4$. Consider the map 
$$
f:\vee_{i=1}^4S^{3k-2}_i \to \vee_{i=1}^4S^k_i
$$
such that $f|S^{3k-2}_i =f_i$ and set $X=C_f$. We define $\ga_m\in H^k(X)$    
the cohomology class which takes the value 1 on the cell $S^k_i$ of $X$ and 0 on 
other cells. For convenience of notation, we set $\ga_5=\ga_1$ and $\ga_6=\ga_2$. 

\begin{lem}\label{basis}
The homology classes $\la \ga_{m}, \ga_{m+1}, \ga_{m+2}\ra$ are 
linearly independent in $H^{3k-1}(X)$.
\end{lem}

\p First, notice all these Massey products are defined and have zero 
indeterminacies. Now, suppose that $\sum_{m=1}^4c_m\la\ga_m, \ga_{m+1}, 
\ga_{m+2}\ra=0$ for some $c_m\in \R$. Consider the space $Z$ as 
in \eqref{z} and the obvious inclusion $j: Z \to X$. Then $j^*\la \ga_{m}, 
\ga_{m+1}, \ga_{m+2}\ra=0$ for $m=2,3,4$ by \corref{zero}, while $j^*\la \ga_1, 
\ga_2, \ga_3\ra\ne 0$ by \theoref{um}. Therefore $c_1=0$. 
Similarly, we can prove that $c_m=0$ for all  $m$.
\qed

\m Now, because of \propref{embed}, $X$ can be regarded as a polyhedron in 
$\R^N$ with $N\ge 4k$. Let $W$ be a regular neighborhood of $X$ in $\R^{N}$ 
and set $M=\partial W$. 

\begin{theorem} If $N\ne 4k,\, 5k-1,\, 6k-2, \,6k-1$ then $H^{3k-1}(M;\Q)$ is 
generated by triple Massey products, while all the cup-products of elements of 
$H^k(M)$ vanish.
\end{theorem}

\p Consider the map
$$
\CD
g: V@>i>> W @>r >> X
\endCD
$$
where $i$ is the inclusion and $r$ is a deformation retraction.
Using the isomorphisms $H^i(W,M)\cong H_{N-i}(X)$ and  $H^i(W)\cong H^i(X)$, and 
the exactness of the sequence
$$
\CD
H^i(W,M) @>>> H^i(W) @>i^*>> H^i(M) @>>> H^{i+1}(W,M).
\endCD
$$
we conclude that 
$H^{2k-1}(M)=0$ and 
$$
g^*: H^{3k-1}(X)\to  H^{3k-1}(M)
$$ 
is an isomorphism. Now, the equality $H^{2k-1}(M)=0$ implies that all the 
Massey products $\la \ga_i, \ga_j, \ga_k \ra$ have zero indeterminacies. 
Furthermore, since $g^*$ is an isomorphism, \lemref{basis} implies that the 
$g^*$-images of the classes $\la \ga_{m}, \ga_{m+1}, \ga_{m+2}\ra$, $m=1,2,3,4$ 
in $M$ form a basis of $H^{3k-1}(M;\Q)$. Finally, the map $i^*: H^k(W) \to 
H^k(M)$ is surjective for $N\ne 4k$, and so the cup-products of elements of 
$H^k(M)$ vanish. 
\qed

\end{document}